\documentclass[12pt]{article}

\title{Sequential adaptive estimators in nonparametric autoregressive models}

\author{Ouerdia Arkoun
\thanks{Laboratoire de Math\'ematiques Rapha\"el Salem,
UMR 6085 CNRS, Universit\'e de Rouen,
Avenue de l'Universit\'e, BP.12,
76801 Saint Etienne du Rouvray (France).\newline
email:  Ouerdia.Arkoun@etu.univ-rouen.fr}
}
\usepackage[round]{natbib}
\usepackage{amssymb}
\usepackage{amsfonts}
\usepackage{amsmath}
\usepackage{amsthm}
\newtheorem{theorem}{Theorem}[section]

\newtheorem{lemma}[theorem]{Lemma}

\newtheorem{remark}[theorem]{Remark}

\newcommand\e{\varepsilon}
\def\bbr{{\mathbb R}}

\newcommand{\wh}{\widehat}

\newcommand\cH{{\cal H}}

\newcommand\cF{{\cal F}}

\newcommand\cN{{\cal N}}
\newcommand\cP{{\cal P}}

\newcommand\cR{{\cal R}}

\def\text#1{\hbox{#1}}
\def\proof{{\noindent \bf Proof. }}

\def\endproof{\mbox{\hfill $\qed$}}

\def\E{{\bf E}}
\def\P{{\bf P}}

\def\C{{\bf C}}

\def\R{{\bf R}}
\def\Chi{{\bf 1}}
\def\d{\mbox{d}}
\def\build #1_#2{\mathrel{\mathop{\kern 0pt #1}\limits_{#2}}} 
\newcommand{\zs}[1]{{\mathchoice{#1}{#1}{\lower.25ex\hbox{$\scriptstyle#1$}}
{\lower0.25ex\hbox{$\scriptscriptstyle#1$}}}}

\numberwithin{equation}{section}

\begin{document}

\maketitle

\begin{abstract}
We construct a sequential adaptive procedure for estimating the autoregressive function at a given point in nonparametric autoregression models with Gaussian noise. We make use of the sequential kernel estimators. The optimal adaptive convergence rate is given as well as the upper bound for the minimax risk.
\end{abstract}

{\bf Key words:}
 Adaptive estimation, kernel estimator, minimax, nonparametric autoregression.\\
\par
{\bf AMS (2000) Subject Classification :}
primary 62G07,62G08; secondary 62G20.

\section{Introduction}

Our problem is the following. Suppose we observe data from the model :
\begin{equation}\label{modeladapt}
y_k=S(x_k)y_{k-1}\, +\xi_k\,,
\quad 1\le k\le n\,, \end{equation}
where $x_k=k/n$ and $(\xi_k)_{k\in\{1,\dots,n\}}$ are random variables independent and identically distributed by standard Gaussian. 

The model \eqref{modeladapt} is a generalization of an autoregressive processes of the first order. In \cite{Da}, the process \eqref{modeladapt} is considered with the function $S,$ having a parametric form. Moreover, the paper of \cite{Da-1} studies spectral properties of the  stationary process \eqref{modeladapt} with the nonparametric function $S$. \cite{Be} considers the model \eqref{modeladapt} with Lipschitz conditions and proposes a recursive estimator. The author establishes the convergence rate for quadratic risk.

This paper deals with a nonparametric estimation of the autoregressive function $S$ at a fixed point $z_0\in]0;1[,$ when the smoothness of $S$ is unknown.
More precisely, we assume that the function $S$ belongs to a H\"{o}lder class with unknown regularity $0 < \beta \leq 1$. Note that for $\beta=1$ this gives the class of Lipschitz functions, considered in \cite{Be}.
The goal of this paper is to find an adaptive minimax convergence rate and to construct an adaptive estimate. 

Many studies is devoted to the minimax convergence rate or an asymptotically efficient estimator in adaptive non sequential setting, i.e. when one or more parameters of the model are assumed to be unknown, in particular, the regularity of the function.
The first result in this direction is obtained in \cite{MR1091202}, where the author proposed an adaptive pointwise estimation method for the Gaussian white noise model. He constructed an adaptive estimation procedure which is minimax for functions from the H\"{o}lder classes with unknown regularity. \cite{MR1867163} modified 
the Lepski{\u\i}'s method for the sequential adaptive estimation for the drift of the diffusion processes.

 In this paper, similarly to \cite{MR1867163}, we apply the Lepski{\u\i} procedure to the model \eqref{modeladapt} based on the sequential kernel estimates. We construct the sequential kernel estimator using the method proposed in \cite{Bor-Kon} for the parametric case. It should be noted that to apply the Lepski{\u\i} procedure the kernel estimators must to have the distribution tail of the Gaussian type. To obtain this property one needs to use the sequential approach. To this end we show some 
modification of the Levy theorem for discrete time and then, using this result, we show that the sequential kernel estimators have the
the same form for the distribution tail as  a Gaussian random variable.  It should be noted that non-sequential kernel estimation does not have the above property 
in the case of the model  \eqref{modeladapt}.  Thus, in this case, the adaptive pointwise estimation is possible only in the sequential framework.
 
  Let we describe now the sequential kernel estimators. For a constant $ H> 0,$ we define $\alpha_\zs{H},$\, $0 \leq \alpha_\zs{H} \leq 1, $ such that

$$\sum^{\tau_\zs{H}-1}_{j=1}\,Q(u_j)\,y^2_{j-1}+ \alpha_\zs{H}\,Q(u_\zs{\tau_\zs{H}})\,y^2_\zs{\tau_\zs{H}-1} = H\,,$$
where the kernel $Q(\cdot)$ is the indicator function on the interval $[-1;1],$ and $\tau_\zs{H}$ is the stopping time defined as follows:
\begin{equation}\label{tauh}
 \tau_H = \inf\{1\leq k\leq n: \sum^k_{j=1}\,Q(u_j)\,y^2_{j-1} \ge H \}.
 \end{equation}
Note that  $$A_\zs{k}=\sum^{k}_{j=1}\,Q(u_j)y_{j-1}^2
\quad\mbox{with}\quad
u_j = \frac{x_j-z_0}{h_n}.
$$
Thus the kernel estimator is written as follows:

\begin{equation}\label{estimateuradapt}
S^*_\zs{H,h_n}(z_0)=\frac{1}{H}\,
\left(\sum^{\tau_\zs{H}-1}_{j=1}\,Q(u_j)\,y_{j-1}\,y_\zs{j}\,+\,\alpha_\zs{H}\,Q(u_\zs{\tau_\zs{H}})\,y_\zs{\tau_\zs{H}-1}\,y_\zs{\tau_\zs{H}}\right)  
\Chi_\zs{(A_\zs{n}\geq H)}.
\end{equation} 
Such an estimator is very convenient to calculate the quantity $\E\,|S^*_\zs{H,h_n}(z_0) - S(z_0)|.$

We describe in detail the statement of the problem in section 2. In section 3 we prove the result of an asymptotic lower bound of adaptive minimax risk. Section 4 is devoted to proving the asymptotic upper bound for the risk of
the kernel estimator \eqref{estimateuradapt}. Section 5 gives the appendix which contains some technical results. Finally, we illustrate the obtained results by numerical examples.

\section{Statement of the problem}

The problem is to estimate the function $S$ at a fixed point $z_\zs{0}\in ]0,1[,$ i.e. the value $S(z_\zs{0})$. For any estimate $\tilde{S}_n = \tilde{S}_\zs{n}(z_0)$ (i.e.  any measurable with respect to the observations $(y_k)_\zs{1\le k\le n}$ function), the risk is defined on the neighborhood $\cH^{(\beta)}(z_\zs{0},K,\e)$ by 
\begin{equation}\label{risqueadapt}
\cR_\zs{n}(\tilde{S}_n)=
\sup_{\beta\in [\beta_\zs{*};\beta^*]}\,\sup_{S\in \cH^{(\beta)}(z_\zs{0},K,\e)}\,N(\beta) \,\E_\zs{S}|\tilde{S}_\zs{n}(z_0)-S(z_0)|\,,
\end{equation}
where $ N(\beta) = \displaystyle\left(\frac{n}{\ln n}\right)^{\beta/(2\beta+1)}$ corresponds to the convergence rate of adaptive estimators on  class $\cH^{(\beta)}(z_\zs{0},K,\e)$ and $\E_\zs{S}$ is the expectation taken with respect to the distribution $\P_\zs{S}$ of the vector $(y_1,...,y_n)$ in \eqref{modeladapt} corresponding to the function $S.$ 

We consider model \eqref{modeladapt} where $S\in \C_\zs{1}([0,1],\bbr)$ is the unknown function. 
 To obtain the stable (uniformly with respect to the function $S$ ) model \eqref{modeladapt}, we assume that for some fixed $0<\e<1,$ the unknown function $S$ belongs to the {\em stability set} 
\begin{equation}\label{3.2}
\Gamma_\zs{\e} = \{S\in \C_\zs{1}(]0,1],\bbr) : \|S\| \le 1-\e \}, 
\end{equation}
where $\|S\| = \sup_{0 < x \leq1}|S(x)|$. Here
$\C_\zs{1}]0,1]$ is the Banach space of continuously differentiable
$]0,1]\to\bbr$ functions.
For fixed  constants $K>0$ and $0<\beta \le 1,$ we define
the corresponding {\em stable local H\"older class} at the point $z_\zs{0}$
as

\begin{equation}\label{3.3}
 \cH^{(\beta)}(z_\zs{0},K,\e) 
= \left\{S \in \Gamma_\zs{\e} \,:\,\Omega^*(z_\zs{0},S) \le K \right\},
\end{equation}
with
$$
\Omega^*(z_\zs{0},S) = \sup_{x\in [0,1]} 
\frac{|S(x) - S(z_\zs{0})|}{|x-z_\zs{0}|^\beta}\,.
$$
The regularity  $\beta \in [\beta_\zs{*};\beta^*]$, is supposed to be unknown, where the interval  $[\beta_\zs{*};\beta^*]$ is known.

%
%
First we give the lower bound for the minimax risk. We show that with the convergence rate $N(\beta)$ the lower bound for the minimax risk is strictly positive.
\begin{theorem}\label{Th.3.2.1}
The risk \eqref{risqueadapt} admits the following lower bound: 
\begin{equation*}
\liminf_{n \to \infty}\,\inf_{\tilde{S}_n}\,\cR_n(\tilde{S}_n) \ge \frac{1}{4}\,,
\end{equation*}
where the infimum is taken over all estimators $\tilde{S_n}.$
\end{theorem}

Now we give the upper bound for the minimax risk of the sequential adaptive estimator defined in \eqref{estimateuradapt}. 
Since $\beta$ is unknown, one can not use this estimator because the bandwidth  $h_n$ depends on $\beta$. That is why we partition the interval $[\beta_\zs{*};\beta^*]$ to follow a procedure of Lepski\u\i. Let us set

\begin{equation}\label{sec:Ad.1}
d_\zs{n} = n/\ln n \quad\mbox{and}\quad  h(\beta)=\left(\frac{1}{d_\zs{n}}\right)^{\frac{1}{2\beta+1}}. 
\end{equation}
We define the grid on the interval $[\beta_\zs{*};\beta^*]$ with the points :
\begin{equation}\label{sec:Ad.2}
\beta_\zs{k}=\beta_\zs{*}+\frac{k}{m}(\beta^*-\beta_\zs{*}),\quad k=0,\dots,m 
\quad\mbox{with}\quad 
m=[\ln d_\zs{n}]+1 
\,.
\end{equation}
We denote $N_\zs{k}$, \,$h_\zs{k}$, \,$S^*_\zs{h}$\, and \, $\omega(h_\zs{j})$\,\, as
$$
N_\zs{k}=N(\beta_\zs{k}),\,\, h_\zs{k}=h(\beta_\zs{k}),\,\, S^*_\zs{h} = S^*_\zs{H,h}\,,
$$
and 
$$
\omega(h_\zs{j})=\max_\zs{0\le k\le j}\,
\left(
|S^*_\zs{h_\zs{j}}
-
S^*_\zs{h_\zs{k}}
|
-
\frac{\lambda}{N_\zs{k+1}}
\right).
$$
We also define the optimal index of the bandwidth as
\begin{equation}\label{sec:Ad.4}
\wh{k}=\inf\left\{0\le j\le m\,:\,\omega(h_\zs{j})\ge \frac{\lambda}{N_\zs{j}}\right\}-1\,.
\end{equation}

We note that $\omega(h_\zs{0})=-\lambda/N_\zs{1}$ and thus $\wh{k}\ge 0$.
The positive parameter, $\lambda,$ is chosen as $\displaystyle \lambda> K+e \sqrt{4+\frac{4}{2\beta_*+1}}.$\\
The adaptive estimator is now defined as
\begin{equation}\label{sec:Ad.5}
\wh{S}_\zs{n}=S^*_\zs{H,\wh{h}}
\quad\mbox{with}\quad
\wh{h}=h_\zs{\wh{k}}\,.
\end{equation}
The following result gives the upper bound for the minimax risk of the sequential adaptive estimator defined above.
\begin{theorem}\label{Th.3.2.2}
 For all $0<\e<1,$ we have
\begin{equation}\label{4.6}
\limsup_{n\to\infty}\,\cR_\zs{n}({\hat{S}}_n) < \infty\,. 
\end{equation}
\end{theorem}

\begin{remark}
Theorem~\ref{Th.3.2.1} gives the lower bound for the adaptive risk, i.e. the convergence rate $N(\beta)$ is best for the adapted risk. Moreover, by Theorem~\ref{Th.3.2.2} the adaptive estimates \eqref{sec:Ad.5} possesses this convergence rate. In this case, this estimates is called optimal in sense of the adaptive risk \eqref{risqueadapt}
\end{remark}

\section{The lower bound}
We show that with this appropriate rate, $N(\beta)$, the lower bound of minimax risk is strictly positive.

{\bf{Proof of Theorem~\ref{Th.3.2.1}}}

To simplify notations, we denote $N(\beta_*) = N_*, \,\, N(\beta^*) = N^*\,\, \mbox{and}\,\, h(\beta_*) = h_*.$\\
We choose $S$ as
$$ S(y) =
\frac{1}{N_*}\,V\left(
\frac{y-z_0}{h_*}\right),$$
where $V$ is a function of $C^{\infty}$ class with compact support $[-1,1]$ such that 
$$ \int_{-1}^{1}\, V^2(u)\, du = \frac{\overline{\beta}}{2} \quad \mbox{with} \quad \overline{\beta} = \frac{\beta^*-\beta_*}{(2\beta^*+1)(2\beta_*+1)}\,, $$
 and satisfying $V(0) = 1 $ and $V(u) = 0$ for $|u|\ge 1.$

It is easy to show that for all real $K,$ large enough,  $S \in \cH^{(\beta_*)}(z_\zs{0},K,\e).$ Note that for all $S,$ the measure $\P_{S}$ is equivalent to the measure $\P_\zs{0}$, where $\P_{0}$ is the distribution of vector $(y_1,\ldots,y_n)$ in \eqref{modeladapt} 
corresponding to function $S_0=0.$
It is also clear that in this case, the density of Radon-Nikodym can be written as

\begin{align*}
\rho_n: &=
  \frac{\d \P_\zs{0}}{\d \P_{S}}(y_1,\dots,y_n)\\
& = \exp\left\{-\frac{1}{2}\sum_{k=1}^n\left(y_k^2 - (y_k-S(x_k)y_\zs{k-1})^2 \right)\right\}\\
& = \exp\left(-\varsigma_n\eta_n-\frac{1}{2}\varsigma_n^2\right),
\end{align*} 

with
$$
\varsigma^2_n=\frac{1}{d_n\,h_*}
\sum^n_{k=1}\,V^2\left(\frac{x_k-z_0}{h_*}\right) y^2_{k-1} 
\quad\mbox{and}\quad
\eta_n=\frac{1}{\sqrt{d_n\,h_*}\,\varsigma_n}\,
\sum^n_{k=1}\,V\left(\frac{x_k-z_0}{h_*}\right)\,y_{k-1}\,\xi_k
\,.
$$
We define 
\begin{equation}\label{tau}
 \tau(S) = 1-S^2(z_0).
\end{equation}
According to Lemma \ref{Le.A.2}, we obtain   
\begin{align*}
\P_{S}-\lim_{n\to\infty}\,\frac{d_n}{n}\,\varsigma^2_n
 &=\P_{S}-\lim_{n\to\infty}\left(\frac{1}{nh_*}\sum^{n}_{k=1}\, V^2\,\left(\frac{x_k-z_0}{h_*}\right) y^2_{k-1}\right)\\
 &=\P_{S}-\lim_{n\to\infty}\frac{1}{\tau(S)}\int_{0}^1 V^2\left(\frac{x-z_0}{h_*}\right) dx\\
 &=\int_{-1}^1 V^2(u) du = \frac{\overline{\beta}}{2} = \varsigma^2_*,
\end{align*}
since $\tau(S) = 1 -\displaystyle\frac{1}{N_*^2}.$

Furthermore, using a central limit theorem for martingales (cf. Lemma \ref{Le.A.6}), it is easy to see that under the measure $\P_S,$ 
$$
 \eta_n\quad \Longrightarrow\quad \cN(0,1)
\quad\mbox{when}\quad n \to \infty
\,.
$$
In fact, we can rewrite $\eta_n$ as follows :

$$\eta_n = \displaystyle \sqrt{\frac{n}{d_n}}\,\frac{\varsigma_*}{\varsigma_n}\,\sum_{k=1}^{n}\,u_{k,n},$$
with $$ u_{k,n}= \frac{1}{\varsigma_*\,\sqrt{n\,h_*}}\,V\left(\frac{x_k-z_0}{h_*}\right)\,y_{k-1}\,\xi_k.$$
Let us consider the first condition of lemma \ref{Le.A.6}. To verify this, it suffices to show that
$$ \E_\zs{S}\,\sum_{k=1}^{n}\,\E_\zs{S}(u^2_{k,n}\Chi_{(|u_{k,n}|>\e)}|\cF_{k-1,n}) \xrightarrow[n\to \infty]{} 0.$$
We have
\begin{align}
\E_\zs{S}\,\sum_{k=1}^{n}\,\E_\zs{S}(u^2_{k,n}\Chi_{(|u_{k,n}|>\e)}|\cF_{k-1,n})
&=\sum_{k=1}^{n}\,\E_\zs{S}\,(u^2_{k,n}\Chi_{(|u_{k,n}|>\e)})\label{cond} \\[3mm]
&= \frac{1}{\varsigma_*^2\,nh_*}\,\sum^{k=k^*}_{k=k_*}\,V^2\,\left(\frac{x_k-z_0}{h_*}\right)\,\E_\zs{S} (y^2_{k-1}\,\xi^2_k\Chi_{(|u_{k,n}|>\e)}),\nonumber
\end{align}
where 
\begin{equation}\label{K}
k_*= [nz_0 - nh_n ] + 1
\quad\mbox{and}\quad
k^*=  [nz_0 + nh_n ]\,,
\end{equation}
with
\begin{align*}
\E_\zs{S} (y^2_{k-1}\,\xi^2_k\Chi_{(|u_{k,n}|>\e)})
&\le \sqrt{\E_\zs{S}\,y^4_{k-1}\,\E_\zs{S}\,\xi^4_k}\,\sqrt{\P_S(|u_{k,n}|>\e)}\\[3mm]
&\le \sqrt{\E_\zs{S}\,y^4_{k-1}\,\E_\zs{S}\xi^4_k}\,\sqrt{\frac{1}{\e^2}\, \E_\zs{S}\,u^2_{k,n}}\\[3mm]
&\le C_1 \sqrt{\frac{\E_\zs{S}\,y^2_{k-1}\,\xi^2_k}{nh_*}} \le \frac{C_2}{\sqrt{nh_*}},
\end{align*}
where $C_1$ and $C_2$ are constants independent of $n.$ So the term in \eqref{cond} is bounded above by
\begin{equation}\label{ineq}
\E_\zs{S}\,\sum_{k=1}^{n}\,\E_\zs{S}(u^2_{k,n}\Chi_{(|u_{k,n}|>\e)}|\cF_{k-1,n}) \leq  \frac{C_3}{nh_*}\,\sum^{k^*}_{k=k_*}\,\frac{1}{\sqrt{nh_*}}\,,
\end{equation}
where $C_3$ is a new constant and as $n\to \infty,$ \eqref{ineq} tends to zero.

\noindent The second condition is easily verified

\begin{align*}
\sum_{k=1}^{n}\,\E_\zs{S}\,(u^2_{k,n}|\cF_{k-1,n}) 
&= \frac{1}{\varsigma^2_*\,n\,h_*}\,\sum_{k=1}^{n}\,V^2\left(\frac{x_k-z_0}{h_*}\right)\,\E(y^2_{k-1}\,\xi_k^2|\cF_{k-1,n}) \\[3mm]
& = \frac{1}{\varsigma^2_*\,n\,h_*}\,\sum_{k=1}^{n}\,V^2\left(\frac{x_k-z_0}{h_*}\right)\,y^2_{k-1}\\
&= \frac{d_n}{n}\,\frac{\varsigma^2_n}{\varsigma^2_*} \xrightarrow[n\to \infty]{\P_S} 1.
\end{align*}

Let us denote $\theta_n = N_*|\tilde{S}_\zs{n}|.$ We have  
\begin{align}
\cR_\zs{n}(\tilde{S}_n) 
& \ge \, \max \left(\E_\zs{S_0}\,N^*|\tilde{S}_\zs{n}|,\E_\zs{S}\,N_*|\tilde{S}_\zs{n}-S(z_0)|\right)\nonumber\\
& = \, \max \left(\E_\zs{S_0}\,\frac{N^*}{N_*}\,|\theta_n|,\E_\zs{S}\,|1-\theta_n|\right)\nonumber\\ 
& \ge \frac{1}{2}\,\E_\zs{S}\,\left(\frac{N^*}{N_*}\,|\theta_n|\,\frac{d\P_0}{d\P_\zs{S}}(y)+ \,|1-\theta_n|\right) \label{3.1} 
\end{align}
%
We set $\displaystyle\gamma_n=\frac{N^*}{N_*}.$ We can rewrite \eqref{3.1} as: 
$$
\cR_\zs{n}(\tilde{S}_n)\ge \frac{1}{2}\,\E_\zs{S}\,(\gamma_n\,\rho_n\,|\theta_n|+ \,|1-\theta_n|). 
$$
Let $B_n = \{\eta_n \le 0\}$ and $C_n = \{\frac{d_n}{n}\,\varsigma^2_n < \overline{\beta} \}.$ Clearly, when $B_n\cap C_n$ is realized, we have
$$\gamma_n\,\rho_n \ge \exp\{ \overline{\beta}\,\ln d_n - \frac{\overline{\beta}}{2} \frac{n}{d_n} \}.$$
The right-hand side of this inequality tends to $\infty$ as $n$ approach $\infty.$
This means that for $n$ sufficiently large,

\begin{align}
\cR_\zs{n}(\tilde{S}_n)
&\ge \frac{1}{2}\,\E_\zs{S}\,\Chi_{{B}_n\cap C_n}(\gamma_n\,\rho_n\,|\theta_n|+ \,|1-\theta_n|)\nonumber\\
&\ge \frac{1}{2}\,\E_\zs{S}\,\Chi_{{B}_n\cap C_n}(|\theta_n|+ \,1-|\theta_n|)\nonumber\\
& =  \frac{1}{2}\,\P_\zs{S}(B_n\cap C_n).\label{3.13}
\end{align}
Since,
$$\P_\zs{S}(B_n\cap C_n) = \P_\zs{S}(B_n) - \P_\zs{S}(B_n\cap C_n^c),$$
$$\P_\zs{S}(B_n\cap C_n^c) \leq \P_\zs{S}(C_n^c) = \P_\zs{S}(\frac{d_n}{n} \varsigma_n \ge \overline{\beta})$$
and $$\frac{d_n}{n}\,\varsigma_n\xrightarrow[n\rightarrow\infty]{\P_S}\frac{\overline{\beta}}{2},$$ hence $$\P_{S}(C_n^c)\xrightarrow[n\rightarrow\infty]{}0.$$\\
As $\P_\zs{S}(B_n) = 1/2,$ we deduce that $\P_\zs{S}(B_n\cap C_n) \xrightarrow[n\rightarrow\infty]{} 1/2.$\\
Passing to the limit as $n\to \infty$ in \eqref{3.13}, we obtain the desired result. \endproof

\section{Sequential adaptive estimation (upper bound)}

{\bf{Proof of Theorem~\ref{Th.3.2.2}}}

We proceed by following a method based on sequential analysis. First, we rewrite the estimation error as follows:
\begin{equation}\label{sec:Up.1}
 S^*_\zs{H,h}(z_0)-S(z_0) = -S(z_0)\,\Chi_\zs{(A_\zs{n}< H)} + B_\zs{H}(h)\,\Chi_\zs{(A_\zs{n}\geq H)}+ \frac{1}{\sqrt{H}}\,\zeta_\zs{H}(h)\,\Chi_\zs{(A_\zs{n}\geq H)}\,,
\end{equation}
where
$$ B_\zs{H}(h) = \frac{1}{H}\,\left( \sum^{\tau_\zs{H}-1}_{j=1}\,Q(u_j)\,(S(x_j)-S(z_0))\,y^2_{j-1}\,+\,\alpha_\zs{H}\,Q(u_\zs{\tau_\zs{H}})\,(S(x_\zs{\tau_\zs{H}})-S(z_0))\,y^2_\zs{\tau_\zs{H}-1}\right) $$ and
 $$ \zeta_\zs{H}(h) = \frac{1}{\sqrt{H}}\,\left( \sum^{\tau_\zs{H}-1}_{j=1}\,Q(u_j)\,y_{j-1}\,\xi_j+\,\alpha_\zs{H}\,Q(u_\zs{\tau_\zs{H}})\,y_\zs{\tau_\zs{H}-1}\,\xi_\zs{{\tau_\zs{H}}}\right).$$ 
Note that the first term in the right-hand side term of  \eqref{sec:Up.1} is studied in Lemma \ref{Le.A.3}. We can show directly that for every   $S\in  \cH^{(\beta)}(z_\zs{0},K,\e)$
\begin{equation}\label{sec:Up.2}
|B_\zs{H}(h)|\le K h^\beta
\end{equation}
 and also, using Lemma \ref{Le.A.5} we have  
\begin{equation}\label{sec:Up.3}
\sup_\zs{n\ge 1}\,
\sup_\zs{h_\zs{*}\le h\le h^*}\,
\E_\zs{S}\,|\zeta_\zs{H}(h)| <\infty\,,
\end{equation}
where $h_\zs{*}=h(\beta_\zs{*})$ and $h^*=h(\beta^*)$.
Now, we choose $H = nh $ and  
$$
\iota=\inf\{k\ge 0\,:\,\beta_\zs{k}\ge \beta\}-1
\,.
$$
This means 
$$
\beta_\zs{\iota}< \beta\le  \beta_\zs{\iota+1}
\quad\mbox{and}\quad
h_\zs{\iota}< h(\beta)\le h_\zs{\iota+1}
\,.
$$
In the sequel, we denote $S^*_\zs{h}(z_0) = S^*_\zs{H,h}(z_0).$ We have now
$$
|S^*_\zs{h_\zs{\iota}}(z_\zs{0})-S(z_\zs{0})|\le  \Chi_\zs{(A_\zs{n}(h_\iota)< nh_\iota)} + K (h(\beta_\zs{\iota}))^\beta+ \frac{1}{\sqrt{nh_\zs{\iota}}}
\,|\zeta_\zs{H}(h_\zs{\iota})|
$$
and
$$
|S^*_\zs{h_\zs{\iota-1}}(z_\zs{0})-S(z_\zs{0})|\le  \Chi_\zs{(A_\zs{n}(h_{\iota-1})< nh_{\iota-1})} + K (h(\beta_\zs{\iota-1}))^\beta +
 \frac{1}{\sqrt{nh_\zs{\iota-1}}} \,|\zeta_\zs{H}(h_\zs{\iota-1})|\,.
$$
Inequality \eqref{sec:Up.3} implies
\begin{equation}\label{sec:Up.4}
\limsup_\zs{n\to\infty}\,
\sup_\zs{\beta_\zs{*}\le \beta\le\beta^*}\,N(\beta)\,
\sup_\zs{S\in  \cH^{(\beta)}(z_\zs{0},K,\e)}\,
\E_\zs{S}\,
\varpi(\iota,z_\zs{0})\,
\,<\infty\,,
\end{equation}
where
$$
\varpi(\iota,z_\zs{0})=|S^*_\zs{h_\zs{\iota-1}}(z_\zs{0})-S(z_\zs{0})|
+|S^*_\zs{h_\zs{\iota}}(z_\zs{0})-S(z_\zs{0})|\,.
$$

Now considering the estimator $\wh{S}_\zs{n},$ one has
\begin{equation}\label{sec:Up.5}
|\wh{S}_\zs{n}(z_\zs{0})-S(z_\zs{0})|\le
I_\zs{1}+I_\zs{2}
+\varpi(\iota,z_\zs{0})\,,
\end{equation}
where
$$
I_\zs{1}=|\wh{S}_n(z_\zs{0})-S(z_\zs{0})|\Chi_\zs{\{\wh{k}\ge \iota+1\}}
\quad\mbox{and}\quad
I_\zs{2}=|\wh{S}_n(z_\zs{0})-S(z_\zs{0})|\Chi_\zs{\{\wh{k}\le \iota-2\}}\,.
$$

We focus now on the left-hand side in this inequality. We have  
\begin{align*}
|\wh{S}_n(z_\zs{0})-S(z_\zs{0})|\Chi_\zs{\{\wh{k}\ge \iota+1\}}\le
|S^*_\zs{\wh{h}}(z_\zs{0})-S^*_\zs{h_\zs{\iota}}(z_\zs{0})|\Chi_\zs{\{\wh{k}\ge \iota+1\}}
+
|
S^*_\zs{h_\zs{\iota}}(z_\zs{0})
-S(z_\zs{0})|\Chi_\zs{\{\wh{k}\ge \iota+1\}}\,.
\end{align*}
Moreover,
\begin{align*}
|S^*_\zs{\wh{h}}(z_\zs{0})&-S^*_\zs{h_\zs{\iota}}(z_\zs{0})|\,\Chi_\zs{\{\wh{k}\ge \iota+1\}}\,
\le
\omega(h_\zs{\wh{k}})\Chi_\zs{\{\wh{k}\ge \iota+1\}}+
\frac{\lambda}{N_\zs{\iota+1}}\\
&\le
\frac{\lambda}{N_\zs{\wh{k}}}\Chi_\zs{\{\wh{k}\ge \iota+1\}}
+
\frac{\lambda}{N_\zs{\iota+1}}
\le
\frac{2\lambda}{N_\zs{\iota+1}}
\le 
\frac{2\lambda}{N(\beta)}
\,.
\end{align*}
This implies directly that
\begin{equation}\label{sec:Up.6}
\limsup_\zs{n\to\infty}\,
\sup_\zs{\beta_\zs{*}\le \beta\le\beta^*}\,N(\beta)\,
\sup_\zs{S\in  \cH^{(\beta)}(z_\zs{0},K,\e)}\,
\E_\zs{S}\,I_\zs{1}\,<\infty\,.
\end{equation}

We establish now a bound for the right-hand side of \eqref{sec:Up.5}:

$$
I_\zs{2}\le 
\left(\Chi_{(A_n(h_{\hat{k}})< nh_{\hat{k}})} +
K (h(\beta_{\hat{k}}))^\beta
+
\frac{1}{\sqrt{nh_{\wh{k}}}}\,\zeta^*
\right)
\Chi_\zs{\{\wh{k}\le \iota-2\}}\,,
$$
where
\begin{equation}\label{4.12}
\zeta^*=\max_\zs{1\le j\le m}\,|\zeta_\zs{H_\zs{j}}(h_\zs{j})|\,.
\end{equation}
Note that
$$
\{\wh{k}\le \iota-2\}=\bigcup^{\iota-1}_\zs{j=1}\,
\left\{\omega(h_\zs{j})\ge \lambda/N_\zs{j}\right\}\,.
$$
Moreover,
\begin{align}
\left\{\omega(h_\zs{j})\ge \lambda/N_\zs{j}\right\}
&=
\bigcup^{j-1}_\zs{l=0}
\left\{
|S^*_\zs{h_\zs{j}}(z_\zs{0})
-
S^*_\zs{h_\zs{l}}(z_\zs{0})
 |
\ge \lambda/N_\zs{j}+\lambda/N_\zs{l+1}\right\} \nonumber \\
&\subseteq \bigcup^{j-1}_\zs{l=0}
\left(
\{
|S^*_\zs{h_\zs{j}}(z_\zs{0})-S(z_\zs{0})|\ge \lambda/N_\zs{j}
\}
\cup
\{
|S^*_\zs{h_\zs{l}}(z_\zs{0})-S(z_\zs{0})|\ge \lambda/N_\zs{l+1}
\}
\right)
\,,\label{union}
\end{align}
also for $j\le \iota-1$,
$$
N_\zs{j}\,(h_\zs{j})^\beta\le \exp\{-\frac{\ln d_\zs{n}}{(2\beta^*+1)m}\}\le 1\,.
$$
For $l\le \iota-1$
$$
N_\zs{l+1}\,(h_\zs{l})^\beta\le \exp\{-\frac{\ln d_\zs{n}}{(2\beta^*+1)m}\}
\le 1,
$$
and
$$
\frac{N_\zs{l}}{N_\zs{l+1}}\ge \exp\{-\frac{\ln d_\zs{n}}{m}\}=e^{-1}\,.
$$
In the first set on the right-hand side \eqref{union}, by Lemma \ref {Le.A.2} we prove that for $n$ sufficiently large and for $\displaystyle \lambda> K+e \sqrt{4+\frac{4}{2\beta_*+1}}$, we have

\begin{align*}
\{|S^*_\zs{h_\zs{j}}(z_\zs{0})-S(z_\zs{0})|\ge \lambda/N_\zs{j}\} 
&\subseteq \left\{ K (h_j)^\beta + \frac{1}{\sqrt{n h_j}} |\zeta_n(h_j)|  \ge \lambda / N_j \right\} \\
&\subseteq \left\{ |\zeta_n(h_j)| \ge  \sqrt{n h_j} \left( \frac{\lambda}{N_j} - K (h_j)^\beta \right)\right\}\,.
\end{align*}
We also have $(1/d_n)^{\beta/(2\beta+1)}\, \sqrt{n h} = \sqrt{n /d_n}$ so that the last inclusion becomes

$$\{|S^*_\zs{h_\zs{j}}(z_\zs{0})-S(z_\zs{0})|\ge \lambda/N_\zs{j}\} \subseteq \left\{|\zeta_n(h_j)| \ge (\lambda - K)\, \sqrt{\frac{n}{d_n}}\right\}.$$

Similarly for the second set on the right-hand side in \eqref{union}, we obtain

$$\{|S^*_\zs{h_\zs{l}}(z_\zs{0})-S(z_\zs{0})|\ge \lambda/N_\zs{l+1}\} \subseteq \left\{|\zeta_n(h_l)| \ge (\lambda- K)/e\, \sqrt{\frac{n}{d_n}}\right\}.$$
Finally,
$$
\{\wh{k}\le \iota-2\}\subseteq
\{\zeta^*\ge \lambda_1\,\sqrt{n/d_\zs{n}} \}\,,
$$
with $\lambda_\zs{1}=(\lambda-K)/e$. So one has
\begin{equation}\label{dernineg}
I_\zs{2}\le \,\Chi_{(A_n(h_{\hat{k}})< nh_{\hat{k}})}\,+\frac{K}{N(\beta)}
+
\frac{1}{\sqrt{nh_\zs{*}}}\,\zeta^*\,
\Chi_\zs{\{\zeta^*\ge \lambda_\zs{1}\,\sqrt{n/d_\zs{n}}\}}.
\end{equation}

Using Lemma \ref {Le.A.2} for $t \ge 2 $, one can easily estimate the first term on the right-hand side of inequality \eqref{dernineg} by

\begin{align*}
\P_\zs{S}(A_n(h_{\hat{k}}) < nh_{\hat{k}})
&= \sum_{l=1}^{m} \P_\zs{S}(A_n(h_l) < nh_{l},\, \hat{k} = l)\\
&\le \sum_{l=1}^{m} \P_\zs{S}(A_n(h_l) < nh_{l})\\
&=  \sum_{l=1}^{m}\,\P_\zs{S}\left(\frac{1}{\tau(S)}\,\int_{-1}^{1}Q(u)du + \Delta_n(Q,h_l) < 1\right)\\
&=   \sum_{l=1}^{m} \,\P_\zs{S} \left( \Delta_n(Q,h_l) < 1-\frac{2}{\tau(S)}\right)\\
&\le  \sum_{l=1}^{m} \,\P_\zs{S} \left( |\Delta_n(Q,h_l)| > 1 \right)\\
& \leq  \sum_{l=1}^{m}\,\E_\zs{S}\,\Delta^{2t}_n(Q,h_l) \le ([\ln d_n]+1)\,C_\zs{1}\,R^{2t}\,(h^*)^{2t\beta}.
\end{align*}
Consider now the last term in  the right-hand side of inequality \eqref{dernineg}. We have
\begin{align*}
 \E_\zs{S}\, \zeta^*\,
\Chi_\zs{\{\zeta^*\ge\lambda_\zs{1}\,\sqrt{\ln n}\}} 
&= \int_{0}^{+\infty}\, \P_\zs{S}(\zeta^*\,\Chi_\zs{\{\zeta^*\ge\lambda_\zs{1}\,\sqrt{\ln n}\}} \ge z ) \, dz \\
&= \int_{0}^{+\infty}\, \P_\zs{S}(\zeta^*\ge z\,, \zeta^*\ge \lambda_\zs{1}\,\sqrt{\ln n})\, dz  \\
&= \lambda_\zs{1}\,\sqrt{\ln n}\,\P_\zs{S}(\zeta^*\ge  \lambda_\zs{1}\,\sqrt{\ln n}) + \int_{\lambda_\zs{1}\,\sqrt{\ln n}}^{+\infty}\,\P_\zs{S}(\zeta^*\ge z)\,dz. 
\end{align*}
Using \eqref{4.12} and Lemma \ref{Le.A.5}, we have
\begin{align*}
\P_\zs{S}(\zeta^*\ge z)
&= \P_\zs{S}(\max_{1\le j\le m}\,|\zeta_n(h_\zs{j})| \ge z)\\
&= \sum_{j=1}^{m}\,\P_\zs{S}(|\zeta_n(h_\zs{j}|) \ge z) \\
&\le 2\, m \, e^{-z^2/8}.
\end{align*}
Then,
\begin{align*}
\E_\zs{S}\,\zeta^*\,
\Chi_\zs{\{\zeta^*\ge\lambda_\zs{1}\,\sqrt{\ln n}\}}
&\le \,2m\,\lambda_\zs{1}\,\sqrt{\ln n}\, e^{-\frac{1}{8}\,\lambda_\zs{1}^2\,\ln n} + 2\,m \int_{\lambda_\zs{1}\,\sqrt{\ln n}}^{+\infty}\,e^{-z^2/8}\,\d z\\
&\le  \,2m\,\lambda_\zs{1}\,\sqrt{\ln n}\, e^{-\frac{1}{8}\,\lambda_\zs{1}^2\,\ln n}+ 2\,m \int_{\lambda_\zs{1}\,\sqrt{\ln n}}^{+\infty}\,z\,e^{- z^2/8}\,\d z \\
& \le \left(\lambda_\zs{1}\,\sqrt{\ln n} +4\right)2m \, n^{-\lambda_\zs{1}^2/8},
\end{align*}
which implies inequality \eqref{4.6}.
\endproof

\section{Appendix}
In this section, we study the properties of stationary processes in the model \eqref{modeladapt}.
\begin{lemma}\label{Le.A.1}
For all $t\in \mathbb{N}^*$ and ~  $0<\e<1,$ the random variables in \eqref{modeladapt} satisfy the following :
\begin{equation}\label{A.1}
r^*=
\sup_\zs{n\ge 1}\,\sup_\zs{0\le k\le n}\,
\sup_\zs{S\in\Gamma_\zs{\e}}\,
\E_\zs{S}\,y^{2t}_\zs{k}\, < \infty.
\end{equation}
\end{lemma}
\proof

Assume that $y_0 = 0.$ Model \eqref{modeladapt} becomes
\begin{equation*}
y_\zs{k}= \sum_{i=1}^{k}\,\prod_{l=i+1}^{k}\, S(x_\zs{l})\,\xi_\zs{i}\,,
\end{equation*} 
 with $S\in\Gamma_\zs{\e}$ and for all $1\le k\le n$,
$$
y^{2t}_\zs{k}\le \left(\sum^k_\zs{j=1}\,(1-\e)^{k-j}\,|\xi_\zs{j}|\right)^{2t}
\,.
$$
Moreover, the H\"older inequality, with  $p = 2t\,,$ gives
\begin{align*}
y^{2t}_\zs{k}
&\le \left(\sum^k_\zs{j=1}(1-\e)^{k-j}\right)^{2t-1}\,\left(\sum^k_\zs{j=1}(1-\e)^{k-j}\,\xi^{2t}_\zs{j}\right) \\
&\le \left(\frac{1}{\e}\right)^{2t-1}\,\left(\sum^k_\zs{j=1}(1-\e)^{k-j}\,\xi^{2t}_\zs{j}\right)\,.
\end{align*}
Thus, it follows that
$$
\E_\zs{S}\,y^{2t}_\zs{k}\,\le\,\frac{(2t)!}{2^t\,t!}\,\left(\frac{1}{\e}\right)^{2t},
$$
and we get the desired result. \endproof\\

Let us introduce the following notation :
$$\Delta_\zs{n}(f,h) = \frac{1}{nh}\,\sum^{n}_{k=1}\,f(u_k)y^2_{k-1}\,-\,\frac{1}{\tau(S)}
\,\int_{-1}^{1}\, f(u)\d u. 
$$
\begin{lemma}\label{Le.A.2}
Let $f$ be a function twice continuously differentiable in $[-1,1]$, 
such that $f(u) =0$ for $|u|> 1$. Then for all $t \in \mathbb{N}^*$,
\begin{equation}\label{A.4} 
 \limsup_{n\to \infty}\,
\sup_\zs{h_*\le h \le h^*}\,\sup_\zs{R>0}
\frac{1}{R^{2t}h^{2t\beta}}\sup_{\|f\|_\zs{1}\le R}\,
 \sup_{S \in \cH^{\beta}(z_0,K,\e)}\,\E_\zs{S}\, 
\Delta^{2t}_\zs{n}(f,h) \le C_\zs{1}\,,
\end{equation}
where
$\|f\|_\zs{1} = \|f\| + \|\dot{f}\|$ \,\, and \,\, $ C_\zs{1} =  2^{4t} K^{2t} (r^*)^2.$
\end{lemma}
\proof
First, write
\begin{equation}\label{A.5} 
\sum^{n}_{k=1}\,f(u_k)y^2_{k-1}=
T_\zs{n}+a_\zs{n}\,,
\end{equation}
where
$$
T_\zs{n}=\sum^{k^*}_{k=k_*}\,f(u_k) y^2_\zs{k}
\quad\mbox{and}\quad
a_\zs{n}=\sum^{k^*}_\zs{k=k_\zs{*}}\,(f(u_\zs{k})-f(u_\zs{k-1}))\,y^2_\zs{k-1}
-f(u_\zs{k^*})\,y^2_\zs{k^*}\,,
$$
 for the integers $k^*$ and $k_\zs{*}$  defined as in \eqref{K}.  Substituting into model \eqref{modeladapt} gives us 
$$
T_\zs{n}= I_\zs{n}(f)+ \sum^{k^*}_{k=k_*}\,f(u_k)S^2(x_k)y^2_{k-1} + M_n\,,
$$
where
 $$
I_\zs{n}(f)=\sum^{k^*}_{k=k_*}\,f(u_k)
\quad\mbox{and}\quad
M_n  =\sum^{k^*}_{k=k_*}\,f(u_k)\,
\left(
2\,S(x_k)\,y_{k-1}\,\xi_k\,+\,\eta_k
\right)
 $$
with $\eta_k = \xi_k^2- 1$.
Noting that
$$
C_n = \sum^{k^*}_{k=k_*}\, (S^2(x_k)-S^2(z_0))\,f(u_k)\,y^2_{k-1}
\quad\mbox{and}\quad
D_n = \sum^{k^*}_{k=k_*}\,f(u_k)(y^2_{k-1}-y^2_k)\,,
$$ 
we obtain 
\begin{equation}\label{A.6} 
\frac{1}{nh}\,
T_\zs{n}
= \frac{1}{\tau(S)}\, \frac{I_\zs{n}(f)}{nh}
+ \frac{1}{\tau(S)} \frac{H_\zs{n}}{nh} 
\end{equation}
with $ H_\zs{n} = M_\zs{n}+C_\zs{n} + S^2(z_0)\,D_\zs{n}$.
Moreover, it is easy to see that
\begin{align*} 
\frac{I_\zs{n}(f)}{nh} &= \int_{-1}^{1}f(t)\d t
 +\sum^{k^*}_{k=k_*}\,\int_{u_{k-1}}^{u_k}f(u_k)\,\d t\,-\, \int_{-1}^{1}f(t)dt \\
& = 
\sum^{k^*}_{k=k_*}\,\int_{u_{k-1}}^{u_k} (f(u_k)- f(t)) dt + \int_{u_{k_{*}-1}}^{u_{k^*}} f(t)dt - \int_{-1}^{1}f(t)dt\,.
\end{align*}
Recall also that $\|f\|+\|\dot{f}\|\le R$. Then 
$$
\left|\frac{1}{nh} \sum^{k^*}_{k=k_*}\,f(u_k)- \int_{-1}^{1}f(t)dt\right| 
 \le \frac{R}{nh}\,.
$$
The definition in \eqref{tau} implies that for any $S\in \Gamma_\zs{\e}$,
\begin{equation}\label{taus}
 \e^2 \leq \tau(S) \leq 1.
\end{equation}

Taking into account \eqref{A.5} and the lower bound for $\tau(S)$  given in \eqref{taus}, we prove that 

\begin{equation}\label{A.7}
\left|\frac{T_\zs{n}}{nh}
-\frac{1}{\tau(S)} \int_{-1}^{1}f(t)dt\right| 
\le 
\frac{1}{\e^2}
\left(
\frac{R}{nh}+\frac{M_\zs{n}}{nh}+
\frac{C_\zs{n}}{nh}+\frac{D_\zs{n}}{nh}
\right)\,.
\end{equation}
We note that $M_\zs{n}$ is the last term of the square integrable martingale  $(G_j)_{k_*\le j\le k^*}$, where
$$ G_j = \sum^{j}_{k=k_*}\,f(u_k)\, \left(2\,S(x_k)\,y_{k-1}\,\xi_k\,+\,\eta_k \right). $$
So, by applying the Burkh\"older inequality, it comes

\begin{align*}
\E_\zs{S}\,\left(\frac{1}{nh}\,M_n\right)^{2t}
& \le \frac{A^{2t}_\zs{2t}}{(n h)^{2t}}\,\E_\zs{S}\,
\left(\sum^{k^*}_{k=k_*}\,f^2(u_k)\,
\left(
2\,S(x_k)\,y_{k-1}\,\xi_k\,+\,\eta_k\right)^2
\right)^{t}\\
& 
\le A^{2t}_\zs{2t}\frac{R^t}{(nh)^{t+1}}\,\E_\zs{S}\,\sum^{k^*}_{k=k_*} \left(
2\,S(x_k)\,y_\zs{k-1}\,\xi_k\,+\,\eta_k\right)^{2t}\\
&\le \frac{R^t}{(nh)^t}\,2^{4t-2}\, A^{2t}_\zs{2t}\,\left(\frac{(2t)!}{2^t t!}\left(2r^*+ \frac{(2t)!}{2^t t!}\right) + 1\right),
\end{align*}
 where $A_\zs{2t} = 18(2t)^{3/2}/(2t-1)^{1/2}$ \,\, and \,\,$r^*$ is given in \eqref{A.1}. 
Since, $|S(x_k)-S(z_0)| \le K|x_k-z_0|^\beta$  
for all $S\in \cH^\beta(z_0,K,\e)$ and after applying the H{\"o}lder inequality for $p=2t$ and $q = 2t/(2t-1),$ 
we obtain
\begin{align*}
 \frac{1}{(n h)^{2t}}\,\E_\zs{S}\, C^{2t}_n
& \le \frac{1}{(n h)^{2t}}\,\left(\sum^{k^*}_{k=k_*}\, |(S^2(x_k)-S^2(z_0))|^{q}\Chi_{|u_k|\le 1}\right)^{2t/q}\,\sum^{k^*}_{k=k_*} f^{2t}(u_k)\,\E_\zs{S}\,y^{4t}_\zs{k-1} \\
& \le 2^{4t}\,R^{2t}\, K^{2t}\,(r^*)^{2}\,h^{2t\beta}\,.
\end{align*}
Now, consider the last term in  the right-hand side of inequality \eqref{A.6}. $D_\zs{n}$ can be written as
$$ 
D_n = \sum^{k^*}_{k=k_*}\,\left((f(u_\zs{k})-f(u_\zs{k-1})\right)\,y^2_\zs{k-1} +f(u_{k_{*}-1})\,y^2_\zs{k_{*}-1} - f(u_\zs{k^*})\,y^2_\zs{k^*}\,.
$$
Since $\|f\|+\|\dot{f}\|\le R$,
we have

$$ \E_\zs{S}\,D^{2t}_\zs{n}
\le\, 2^{4t-2}\,R^{2t}  
\E_\zs{S}\,
\left(\frac{1}{n h}
\,
\sum^{k^*}_{k=k_*}\,y^{4t}_\zs{k-1}\,
+ y^{4t}_\zs{k^*}+ y^{4t}_\zs{k_{*}-1})
\right)
 \le 2^{4t}\,R^{2t}\,(r^*)^2\,.
$$ 
 Similarly we find a bound for the second term in  the right-hand side of the expression \eqref{A.4}.
Hence we have Lemma ~\ref{Le.A.2}. 
\endproof

\begin{lemma}\label{Le.A.3}
For any $t\ge 1$, the stopping time $\tau_\zs{H}$ defined in \eqref{tauh} satisfies the following property: for $H = nh$
$$\P_\zs{S} (\tau_H > n) \le C_\zs{1}(Rh)^{2t\beta }\,,$$
where $C_\zs{1}$ is defined in \eqref{A.4}.
\end{lemma}

\proof

Taking into account that $\tau(S)\le 1$, we obtain 
\begin{align*}
\P_\zs{S} (\tau_H > n) 
&= \P_\zs{S}(\frac{1}{nh}\sum_{k=1}^{n} Q(u_k)\,y^2_{k-1} < \frac{H}{nh})\\
&= \P_\zs{S}\left(\frac{1}{\tau(S)}\,\int_{-1}^{1}Q(u)du + \Delta_n(Q,h) < 1\right)\\
&=  \P_\zs{S} \left( \Delta_n(Q,h) < 1-\frac{2}{\tau(S)}\right)\\
&\le \P_\zs{S} \left( |\Delta_n(Q,h)| > 1 \right) \leq \E_\zs{S}\,\Delta^{2t}_n(Q,h) \le C_\zs{1}\,R^{2t}\,h^{2t\beta}\,.
\end{align*}
This last inequality comes from Lemma \ref{Le.A.2}
\endproof\\

To prove Lemma \ref{Le.A.5}, we need the following lemma proved in \cite{MR0488267} p.234-235.
\begin{lemma}\label{Le.A.4}
Let the Wiener process  $W = (W_t,\cF_t),$\,$t\ge 0,$ be given on a probability space and let there also be given the random process $f = (f_t,\cF_t),$\, $t\ge 0,$ such that :\\
 $$ (1) \quad  \quad \quad  P\left(\int_{0}^{T}\, f^2_t dt < \infty \right) = 1, \,\,\,\, 0<T<\infty\,,$$
 $$ (2) \quad \quad  \quad \quad \quad \quad \quad \quad \quad   P\left(\int_{0}^{\infty}\, f^2_t dt = \infty \right) = 1. $$
Then the random process $z = (z_\zs{s},\Gamma_\zs{s}),\, s\ge 0,$ with $z_\zs{s} = \int_{0}^{\tau_\zs{s}}\, f_t\, dW_t,\,\Gamma_\zs{s} = \cF_{\tau_\zs{s}},$ where $\tau_\zs{s}= \inf(t:\int_{0}^{t}\, f^2_u du > s),$ is a Wiener process and with probability one 

$$\lim_{t\to \infty}\, \frac{\int_{0}^{t}\, f_u dW_u}{\int_{0}^{t}\, f^2_u du} = 0.$$
\end{lemma}

\begin{lemma}\label{Le.A.5}
For all $z\ge 2$ and $H>0,$ one has   
\begin{equation}
 \P_\zs{S}(|\zeta_H(h)|>z)\le 2\, e^{-z^2/8}.
\end{equation}
\end{lemma}

\proof
The Brownian motion, $(W_t)_{t\ge 0},$ is a stochastic process whose disjoint increments are independent as $W_{t+s} - W_t$, follows a Gaussian distribution with zero mean and variance $s$. So in our case, we can write
 $$\xi_k = W_k-W_{k-1}\,\,  \sim \cN(0,1).$$
We recall that
$$  \zeta_\zs{H} = \frac{1}{\sqrt{H}}\left( \sum^{\tau_\zs{H}-1}_{j=1}\,Q(u_j)\,y_{j-1}\,\xi_j+\,\alpha_\zs{H}\,Q(u_\zs{\tau_\zs{H}})\,y_\zs{\tau_\zs{H}-1}\,\xi_\zs{{\tau_\zs{H}}} \right)\Chi_\zs{(A_\zs{n}\ge H)}. $$
Thus,
$$\P_\zs{S}(|\zeta_\zs{H}| > z\,\, \Chi_\zs{(A_\zs{n}\ge H)}) = \P_\zs{S}(|\zeta_\zs{H}| > z , \Chi_\zs{(A_\zs{n}\ge H)}) = \P_\zs{S}(|\tilde{\zeta}_\zs{H}| > z , \Chi_\zs{(A_\zs{n}\ge H)}),$$
 where 
 $$ \tilde{\zeta}_\zs{H}(h) = \frac{1}{\sqrt{H}}\left( \sum^{\tilde{\tau}_\zs{H}-1}_{k=1}\,\delta_k \, \xi_k +\,\tilde{\alpha}_\zs{\tilde{\tau}_\zs{H}}\,\delta_{\tilde{\tau}_\zs{H}} \,\xi_\zs{\tilde{\tau}_\zs{H}} \right)$$
 and $$ \sum^{\tilde{\tau}_\zs{H}-1}_{k=1}\,\delta^2_\zs{k} \, + \,\tilde{\alpha}_\zs{\tilde{\tau}_\zs{H}}\,\delta^2_{\tilde{\tau}_\zs{H}} = H, $$
 with \,$\delta_k = Q(u_k)\,y_{k-1}\,\Chi_\zs{(k \le k^*)} + \Chi_\zs{(k > k^*)} $ \,\,and
 $$\tilde{\tau}_\zs{H} = \inf\{ k\ge 1 : \sum^{k}_{j=1}\, \delta_\zs{j}^2 \ge H \}. $$
One can see that
$$\P_\zs{S}(|\tilde{\zeta}_\zs{H}(h)| > z) = \P_\zs{S}\left(\frac{1}{\sqrt{H}}\,\left|\int_{0}^{\tilde{\tau}_\zs{H}}\, f_\zs{t}\, \d W_\zs{t}\right| > z \right),$$
where $$f_\zs{t} = \sum^{\infty}_{j=1}\, \delta'_\zs{j}\,\Chi_\zs{[j-1,j]}(t) $$
with
\begin{equation}
\delta'_\zs{j} =
\left\{
\begin{array} {ll}
 \delta_j & j< \tilde{\tau}_\zs{H}\\
\tilde{\alpha}_\zs{\tilde{\tau}_\zs{H}}\,\delta_\zs{\tilde{\tau}_\zs{H}}  &  j= \tilde{\tau}_\zs{H}\\
0   &    j > \tilde{\tau}_\zs{H}. 
\end{array}
\right.
\end{equation}
Indeed, 
\begin{align*}
\int_{0}^{\tilde{\tau}_\zs{H}}\, f_\zs{t}\, \d w_\zs{t} 
= \sum^{\tilde{\tau}_\zs{H}}_{j=1}\,\int_{j-1}^{j}\, f_\zs{t}\, \d w_\zs{t}
&=  \sum^{\tilde{\tau}_\zs{H}}_{j=1}\, \delta'_\zs{j}\,[w_\zs{j} - w_\zs{j-1}]\\
&=  \sum^{\tilde{\tau}_\zs{H}}_{j=1}\, \delta'_\zs{j}\,\xi_\zs{j} = \sum^{\tilde{\tau}_\zs{H}-1}_{j=1}\,\delta_\zs{j} \, \xi_\zs{j} +\,\tilde{\alpha}_\zs{\tilde{\tau}_\zs{H}}\,\delta_{\tilde{\tau}_\zs{H}} \,\xi_\zs{\tilde{\tau}_\zs{H}}. 
\end{align*}
We set
$$g_\zs{t} = \sum^{\infty}_{j=1}\, \delta''_\zs{j}\,\Chi_\zs{[j-1,j]}(t) $$
with
\begin{equation}
\delta''_\zs{j} =
\left\{
\begin{array} {ll}
 \delta_j & j< \tilde{\tau}_\zs{H}\\
\sqrt{\tilde{\alpha}_\zs{\tilde{\tau}_\zs{H}}}\,\delta_\zs{\tilde{\tau}_\zs{H}}  &  j= \tilde{\tau}_\zs{H}\\
0   &    j > \tilde{\tau}_\zs{H}\,,
\end{array}
\right.
\end{equation}
which yields,
\begin{align*}
\int_{0}^{\tilde{\tau}_\zs{H}}\, g^2_\zs{t}\, \d t 
&= \sum^{\tilde{\tau}_\zs{H}}_{j=1}\,\int_{j-1}^{j}\, g^2_\zs{t}\, \d t \\
& = \sum^{\tilde{\tau}_\zs{H}-1}_{j=1}\,\delta^2_\zs{j} \, + \,\tilde{\alpha}_\zs{\tilde{\tau}_\zs{H}}\,\delta^2_{\tilde{\tau}_\zs{H}} = H.
\end{align*}
By lemma \ref{Le.A.4}, we obtain
 $$\eta = \frac{1}{\sqrt{H}} \int_{0}^{\tilde{\tau}_\zs{H}}\, g_\zs{t}\, \d W_\zs{t} \sim \cN(0,1).$$
Then,
\begin{eqnarray}
& & \P_\zs{S}\left(\frac{1}{\sqrt{H}}\,\left|\int_{0}^{\tilde{\tau}_\zs{H}}\, f_\zs{t}\, \d W_\zs{t}\right| > z \right) \nonumber\\[3mm]
& \leq & \P_\zs{S}\left(\frac{1}{\sqrt{H}}\,\left|\int_{0}^{\tilde{\tau}_\zs{H}}\, g_\zs{t}\, \d W_\zs{t}\right| > \frac{z}{2} \right) + \P_\zs{S}\left(\frac{1}{\sqrt{H}}\,\left|\int_{0}^{\tilde{\tau}_\zs{H}}\, (f_\zs{t} - g_\zs{t})\, \d W_\zs{t}\right| > \frac{z}{2} \right) \nonumber \\ [3mm]
& \le & \P_\zs{S}\left(|\eta| > \frac{z}{2} \right) + \P_\zs{S}\left(\frac{1}{\sqrt{H}}\,|\sqrt{\alpha_\zs{\tilde{\tau}_\zs{H}}} - \alpha_\zs{\tilde{\tau}_\zs{H}}| |\delta_\zs{\tilde{\tau}_\zs{H}}\,\xi_\zs{\tilde{\tau}_\zs{H}}| > \frac{z}{2} \right) \nonumber \\ [3mm] 
& \le & \P_\zs{S}\left(|\eta| > \frac{z}{2} \right) + \P_\zs{S}\left(\frac{1}{\sqrt{H}}\,\sqrt{\alpha_\zs{\tilde{\tau}_\zs{H}}} \,|\delta_\zs{\tilde{\tau}_\zs{H}}|\,|\xi_\zs{\tilde{\tau}_\zs{H}}| > \frac{z}{2} \right)\nonumber  \\[3mm]  
& = & \P_\zs{S}\left(|\eta| > \frac{z}{2} \right) + \P_\zs{S}\left(\frac{1}{H}\,\alpha_\zs{\tilde{\tau}_\zs{H}} \,\delta^2_\zs{\tilde{\tau}_\zs{H}}\,\xi^2_\zs{\tilde{\tau}_\zs{H}} > \frac{z^2}{4} \right)\nonumber  \\ [3mm]
& \le & \P_\zs{S}\left(|\eta| > \frac{z}{2} \right) + \P_\zs{S}\left(\xi^2_\zs{\tilde{\tau}_\zs{H}} > \frac{z^2}{4} \right).\label{A.12}
\end{eqnarray}
As $\eta$ is a standard Gaussian random variable, we can write for all $z\ge 2$, 
\begin{align*} 
\P_\zs{S}\left(|\eta| > \frac{z}{2} \right)
& = \sqrt{\frac{2}{\pi}}\int_{z/2}^{+\infty} e^{-t^2/2} dt \\[3mm]
&\leq  \sqrt{\frac{2}{\pi}}\int_{z/2}^{+\infty} t\,e^{-t^2/2} dt = \sqrt{\frac{2}{\pi}}\, e^{-z^2/8}.
\end{align*}
We can also express the second term on the right-hand side of \eqref{A.12} as  
\begin{align*}
\P_\zs{S}\left(\xi^2_\zs{\tilde{\tau}_\zs{H}} > \frac{z^2}{4} \right)
& = \sum_{l =1}^{+ \infty}\,\P_\zs{S}\left(\xi^2_\zs{l} > \frac{z^2}{4}\,,\,\tilde{\tau}_\zs{H} = l \right)\\[3mm]
&= \sum_{l =1}^{+ \infty}\,\P_\zs{S}\left(\xi^2_\zs{l} > \frac{z^2}{4} \,,\,\sum_{j=1}^{l -1}\,\delta^2_\zs{j}<H\,,\, \sum_{j=1}^{l }\,\delta^2_\zs{j}\ge H  \right)\\[3mm]
&= \sum_{l =1}^{+ \infty}\,\P_\zs{S}\left(|\xi_\zs{l}| > \frac{z}{2}\right)\,\P_\zs{S}(\tilde{\tau}_\zs{H} = l )\\[3mm]
&\le \sqrt{\frac{2}{\pi}}\, e^{-z^2/8}\,\sum_{l =1}^{+ \infty}\,\P_\zs{S}(\tilde{\tau}_\zs{H} = l) = \sqrt{\frac{2}{\pi}}\, e^{-z^2/8}.
\end{align*}
So for any $z \ge 2,$ \eqref{A.12} implies

$$\P_\zs{S}\left(\frac{1}{\sqrt{nh}}\,\left|\int_{0}^{\tilde{\tau}_\zs{H}}\, f_\zs{t}\, \d w_\zs{t}\right| > z \right) \le 2\,\sqrt{\frac{2}{\pi}} \,e^{-z^2/8}.$$
\endproof

\begin{lemma}\label{Le.A.6}\citep[][pp. 80-82]{He}~\\
Let $(u_{k,n})_{1\le k\le n}$ be a "martingale difference" defined on a probability space $(\Omega,\cF, \cP)$ and filtrations $\{\cF_{k,n}, k\in \mathbb{N}\}$ of $\cF,$ $n\in \mathbb{N}^*$ such that $u_{k,n}$ is $\cF_{k,n}$-measurable. Assume that the following two conditions are satisfied:

 $$ (1)  \quad  \quad  \quad \sum_{k=1}^{n}\, \E(u^2_{k,n}\Chi_{(|u_{k,n}|>\e)}|\cF_{k-1,n}) \xrightarrow[n\to \infty]{\P} 0, \quad \mbox{for all} \quad \e > 0, $$

 $$ (2)  \quad  \quad  \quad  \quad \quad  \quad  \quad  \quad \quad  \quad  \quad  \quad  \quad  \quad \sum_{k=1}^{n} \,\E(u^2_{k,n}|\cF_{k-1,n}) \xrightarrow[n\to \infty]{\P} 1.$$

Then, $$\sum_{k=1}^{n} u_{k,n} \Longrightarrow \cN(0,1).$$
\end{lemma}
\section{Numerical simulations }
We illustrate the obtained results by the following simulation which is established using Scilab.

The purpose is to estimate, at a given point $z_0,$ the function $S$ defined over $[0;1]$ by $S(x) = |x-z_0|^\beta.$
We check that such a function belongs to $\cH^{(\beta)}(z_\zs{0},K,\e)$ when $K \ge 1.$
The values of  $z_0$ and $\beta$ are arbitrary, which permit the user to name his choice. As an example, take $z_0=1/\sqrt{2}.$ Then $\beta_*=0.6$ is a lower regularity value and $\beta^*=0.8$ is the higher regularity value.

We simulated $n$ data for the function $S(x)=|x-z_0|^{\beta}$ for $\beta = 0.7.$ We obtained an estimation in constructing the estimator $\hat{S}_n$  defined in \eqref{sec:Ad.5} with the procedure of Lepski\u\i~ which gives us the optimal window for the index $\hat{k}$ defined in \eqref{sec:Ad.4}.

Numerical results approximate the asymptotic risk of a sequential estimator defined in \eqref{sec:Ad.5} used due to the calculation of an expectation (it performs an average for  $M=15000$ simulations) and the finite number of observations $n$. Here we calculate for the sequential estimator the quantity $\R_n = \displaystyle\frac{1}{M}\,\sum_{k=1}^{M}|\hat{S}_n^{(k)}(z_0) - S(z_0)|.$

By varying the number of observations $n,$ we obtain different risks listed in the following table:

\begin{center}
\begin{tabular}{|c|c|c|c|c|}
  \hline
  $n$ &        100 &  1000 &  5000& 10000 \\
  \hline
 $\R_n$ & 0.284 & 0.154 & 0.101 & 0.087 \\
  \hline
  
\end{tabular}
\end{center}
 
 When taking $\beta = \beta^* = 1$, we obtain
\begin{center}
\begin{tabular}{|c|c|c|c|c|}
  \hline
  $n$ &        100 &  1000 &  5000& 10000 \\
  \hline
 $\R_n$ & 0.201 & 0.097 & 0.058 & 0.047 \\
  \hline
  
\end{tabular}
\end{center}

\bibliographystyle{plain}


\end{document}